\documentclass[letterpaper]{article}
\usepackage[letterpaper,margin=1in]{geometry}

\usepackage{amsmath,amsfonts,amssymb,amsthm}

\usepackage[colorlinks,citecolor=red]{hyperref}

\newtheorem{definition}{Definition}[section]
\newtheorem{theorem}{Theorem}[section]
\newtheorem{lemma}{Lemma}[section]
\newtheorem{remark}{Remark}[section]
\newtheorem{corollary}{Corollary}[section]
\newtheorem{proposition}{Proposition}[section]

\newcommand{\RR}{\mathbb{R}}
\newcommand{\NN}{\mathbb{N}}

\newcommand{\norm}[1]{\left|#1\right|}

\DeclareMathOperator*{\esssup}{ess\,sup}

\title{A non-Gaussian Hardy-type Equation in Fractional Time}

\author{Soveny Solís and Vicente Vergara}
\date{}

\begin{document}

\maketitle

\abstract{A non-Gaussian Hardy equation is studied with a non-linearity of Osgood-type growth. A fractional derivative in time is incorporated for the first time in an research of this type. Existence of local and global solutions are established by combining properties of the fundamental solutions together with the parameters of the non-Gaussian process, leading to optimal asymptotic estimates. Additional properties of the fundamental solutions and instantaneous blow-up results are found. The Banach contraction mapping principle is particularly exploited. It is also defined a critical exponent for existence and non-existence of solutions together with a judicious choice of the initial data.}


\begin{center}
\end{center}




\section{Introduction}

Let $\alpha\in (0,1)$, $\beta\in (0,2)$, $p>1$ and $\gamma >0$ fixed parameters. We investigate the non-Gaussian parabolic Hardy-type equation in the form
\begin{equation}\label{FE:0}
	\partial_t^{\alpha} (u - u_0) + \Psi_{\beta}(-i\nabla) u = |x|^{-\gamma}|u|^{p-1}u,\quad t>0,\ x\in \RR^d,
\end{equation}
where $\partial_t^{\alpha}$ denotes the Riemann-Liouville fractional derivative of order $\alpha$. This derivative is defined as:
\[
\partial_t^{\alpha} v = \frac{d}{dt}\int_0^t g_{1-\alpha}(t-s)v(s)ds = \frac{d}{dt} (g_{1-\alpha} \ast v)(t),
\]
with $g_{\rho}(t)=\frac{1}{\Gamma(\rho)}t^{\rho-1}$, where $\Gamma(\cdot)$ represents the gamma function.

The operator $\Psi_{\beta}(-i\nabla)$ is a pseudo-differential operator of order $\beta$, with symbol $\psi_{\beta}(\xi) = |\xi|^{\beta} \omega_{\nu}(\xi/|\xi|)$ for $\xi\in \RR^d$. Here, $\omega_{\nu}$ is a continuous function on the surface of the unit sphere $\mathbb{S}^{d-1}$, and $\nu$ denotes a spectral measure on $\mathbb{S}^{d-1}$. Specifically:
\begin{equation}\label{omega:nu}
	\omega_{\nu}(\theta) :=\int_{\mathbb{S}^{d-1}} |\theta \cdot \eta|^{\beta}\nu(d\eta), \quad \theta\in \mathbb{S}^{d-1}.
\end{equation}
For further details regarding the operator $\Psi_{\beta}(-i\nabla)$, see \cite{Kolo09}.

The evolution problem \eqref{FE:0} is considered throughout this work under the following general hypotheses for the initial data $u_0$ and the spectral measure $\nu$:
\begin{itemize}
	\item[(H1)] The initial data $u_0$ belongs to $L_q(\RR^d)$, with $1 \leq q < \infty$.
	\item[(H2)] The spectral measure $\nu$ has a strictly positive density, such that the function $\omega_{\nu}$ in \eqref{omega:nu} is strictly positive and differentiable $(d + 1 + [\beta])$ times on $\mathbb{S}^{d-1}$.
\end{itemize}

Condition $(H2)$ guarantees a two-sided estimate for the fundamental solution of the linear part of \eqref{FE:0}, as seen in \eqref{z:doble:est} and \eqref{y:doble:est} below. We exploit these estimates to derive our main results.

We can express the solution to \eqref{FE:0} using the \textit{variation formula for Volterra equations} (cf. \cite{Pr93}) as follows:
\[
u(t,x) := S(t)u_0(x) + \int_{0}^{t} R(t-s)\left(|\cdot|^{-\gamma}|u|^{p-1}u(s)\right)(x) ds,
\]
where, for each $t>0$, $S(t)$ and $R(t)$ are linear and bounded operators on $L_p(\RR^d)$ of convolution type, defined as:
\begin{equation}
\label{operatosRS}
S(t)v(x) := \int_{\RR^d} Z(t,x-y)v(y)dy, \quad \text{and} \quad R(t)v(x):= \int_{\RR^d} Y(t,x-y) v(y)dy.
\end{equation}
When dealing with positive solutions, it is noteworthy that in the limiting cases as $\alpha \to 1$ and $\beta\to 2$, with $\omega_{\nu}\equiv 1$, equation \eqref{FE:0} reduces to the classical parabolic Hardy-type equation:
\begin{equation}\label{FE:0:1:2}
	\partial_t u - \Delta u = |x|^{-\gamma}u^p,\quad t>0,\ x\in \RR^d, \text{ and } u(t,x)|_{t=0} = u_0(x), \, x\in \RR^d.
\end{equation}
In this case, the operators $S$ and $R$ coincide with heat semigroup, i.e.,
\[ 
S(t)v = R(t)v = \frac{1}{(4\pi t)^{d/2}}\int_{\mathbb{R}^d} e^{-|x-y|^2/4t}v(y) \, dy, 
\]
since the kernels $Z$ and $Y$ match the Gaussian heat kernel $p_t(x) := \frac{1}{(4\pi t)^{d/2}}e^{-|x|^2/4t}$. However, despite this alignment and the use of the Gaussian kernel and the properties of the heat semigroup, extensively investigated in the literature, the study of the properties of the solution to \eqref{FE:0:1:2} presents non-trivial challenges, as demonstrated in \cite{HiTa21}. In our case, where $Z$ and $Y$ are generally distinct and not explicitly known, we aim to demonstrate properties of the solution to \eqref{FE:0} similar to those obtained for \eqref{FE:0:1:2}, at least partially. To provide insight into the types of kernels $Z$ and $Y$, we will briefly summarize the state-of-the-art in this subject. 

The fundamental solution, denoted by $Z(t,x)$, for \eqref{FE:0} arises when setting $u_0$ as the Dirac measure and the right-hand side as 0 in \eqref{FE:0}. For $\Psi_{\beta}(-i\nabla)=-\Delta$ with $\beta=2$ and $\omega_{\nu}\equiv 1$, and $\alpha\in (0,1)$, it is known (see, e.g., \cite{ScWy89, Koch90}) that
\[
Z(t,x)=\pi^{-\frac{d}{2}}t^{\alpha-1}|x|^{-d} H^{20}_{12}\left(\frac{1}{4}|x|^2t^{-\alpha}\big|^{(\alpha,\alpha)}_{(d/2,1), (1,1)}\right),\quad t>0,\,x\in \RR^d\setminus\{0\},
\]
where $H$ represents the Fox $H$-function \cite{KiSa04, KiST06}. Nevertheless, this representation of $Z$ is not conducive to derive direct estimates due to the complexity of the $H$-function. By utilizing the analytic and asymptotic properties of $H$, \cite{EiKo04} (and \cite{Koch90}) derived sharp estimates for $Z$. Alternatively, in \cite{KeSVZ16}, $Z(t,x)$ was obtained using the subordination principle for abstract Volterra equations with completely positive kernels, which is detailed in \cite[Chapter 4]{Pr93} and \cite{ClNo79} (also \cite{PoVe18} and references therein). Specifically, $Z(t,x)$ is derived from the heat kernel as follows:
\begin{equation*}
	Z(t,x) = -\int_0^{\infty} p_{s}(x)\omega(t,ds),\quad t>0, x\in \RR^d,
\end{equation*}
where $-\omega(t,ds)$ is a probability measure on $\RR_+$ for each $t>0$.

In \cite{KeSZ17} (also in \cite{KiLi16}), $Z$ is considered when $\Psi_{\beta}( -i\nabla)=(-\Delta)^{\beta/2}$ with $\beta\in (0,2)$ and $\alpha\in (0,1)$. The method from \cite{EiKo04} is used to construct and estimate the corresponding kernel $Z$. However, \cite[Section 2]{JoKo19} and \cite[Section 8.2]{Kolo09} provide crucial developments on this topic. It is shown that the linear Cauchy problem \eqref{FE:0} admits a fundamental solution $Z$, given by
\begin{equation}\label{Z:G}
	Z(t,x)=\dfrac{1}{\alpha}\displaystyle\int_0^\infty G(t^\alpha s,x)s^{-1-\frac{1}{\alpha}}G_\alpha(1,s^{-\frac{1}{\alpha}})ds,
\end{equation}
where $G$ is the Green function solving $\partial_t v(t,x)+\Psi_{\beta}(-i\nabla)v(t,x)=0$ with $v(0,x)=\delta_0(x)$, and $G_{\alpha}$ solves $\partial_t v(t,s)+\dfrac{d^\alpha}{ds^\alpha}v(t,s)=0$ with $v(0,s)=\delta_0(s)$, with $\dfrac{d^\alpha}{ds^\alpha}f(s):=\dfrac{1}{\Gamma(-\alpha)}\int_0^\infty \dfrac{f(s-\tau)-f(s)}{\tau^{1+\alpha}}d\tau$ as given in \cite{Kolo19}. By combining \eqref{Z:G} with Aronson's approach \cite{Aron67}, authors in \cite{JoKo19} derive double-sided estimates for $Z(t,x)$ under condition (H2) on $\nu$. This representation of $Z$ offers several advantages over previous ones, being applicable to a wider range of operators $\Psi_{\beta}(-i\nabla)$, including fractional Laplacians of order $\beta\in(0,2)$. Moreover, it is expressed in terms of the more intuitive Green function $G$, and is more amenable to analysis, facilitating explicit estimates for $Z$ in certain cases. More precisely, for a fixed $T > 0$ and $(t, x, y) \in (0, T] \times \RR^d \times \RR^d$, the following two-sided estimates for $Z(t,x-y)$ hold. Considering $\Omega:=|x-y|^{\beta} t^{-\alpha}\leq 1$ we have
\begin{equation}\label{z:doble:est}
Z(t,x-y) \asymp 
\begin{cases}
	t^{-\frac{d\alpha}{\beta}} & d<\beta,\\
	t^{-\alpha} (|\log \Omega| + 1) & d = \beta,\\
	t^{-\frac{d\alpha}{\beta}} \Omega^{1-\frac{d}{\beta}} & d > \beta,
\end{cases}
\end{equation}
and for $\Omega\geq 1$
\[
Z(t,x-y) \asymp C t^{-\frac{d\alpha}{\beta}} \Omega^{-1-\frac{d}{\beta}}.
\]
The notation $f(x)\asymp g(x)$ in $D$ above means that there exists constants $C, c > 0$ such that $f$ satisfies the following two-sided estimate, $cg(x) \leq f(x) \leq Cg(x)$, for all $x \in D$, for some region $D$.

Regarding the kernel $Y$, it is defined as the unique solution of the equation
\begin{equation}\label{Z:Y}
	Z(t,x) = (g_{1-\alpha}\ast Y(\cdot,x))(t),
\end{equation}
which also enjoys of two-sided estimates similar to \eqref{z:doble:est} as follows: for $\Omega\leq 1$,
\begin{equation}\label{y:doble:est}
	Y(t,x-y) \asymp 
	\begin{cases}
		t^{-\frac{d\alpha}{\beta}+\alpha-1} & d<2\beta,\\
		t^{-\alpha-1} (|\log \Omega| + 1) & d = 2\beta,\\
		t^{-\frac{d\alpha}{\beta}+\alpha-1} \Omega^{2-\frac{d}{\beta}} & d > 2\beta,
	\end{cases}
\end{equation}
and for $\Omega\geq 1$ 
\[
Y(t,x-y) \asymp t^{-\frac{d\alpha}{\beta}+\alpha-1} \Omega^{-1-\frac{d}{\beta}},
\]
see \cite[Proposition 2.2]{SoVe22} and \cite[Lemma 2.15]{SoVe22} for the derivation of \eqref{Z:Y}.

Modeling of diffusive processes has been a topic of interest in different areas of science. This is why pseudo-differential equations have gained ground in recent years, in particular, non-Gaussian diffusive equations. Together with fractional derivatives in time, these new models constitute the theoretical base for the study of non-local phenomena in physics, finance, biology, and other fields, see e.g. \cite{ADK23}, \cite{AwMe20}, \cite{Bel21}, \cite{KiLi16}, \cite{LW21} and \cite{SaTo19}.

In this work, one of our goals is to study the so-called \textit{critical exponent} for the Hardy equation \eqref{FE:0} in a non-Gaussian context. The value $q_c$ given in \eqref{Cexponent} below, is considered the critical exponent in the sense that initial data is chosen depending on the $\frac{d}{q_c}$ radio, which can lead to a well-posed or ill-posed problem. 

When $\gamma=0$, other studies for \eqref{FE:0} have been carried out. For instance, existence of global solutions can be found in  \cite{SoVe22} and \cite{ZS15}, blow-up in finite time in \cite{SoVe23} and non-existence of solutions in \cite{SoVe25}.

The case $\gamma >0$ has been investigated in other particular contexts. Thus, we see in the literature many efforts dedicated to the Hardy equation for $\beta =2$ and $\alpha =1$ and $\omega_{\nu}(\theta) \equiv 1$, see e.g., \cite{STW17}, \cite{HiTa21} and references therein. The case of fractional Laplacian in  \eqref{FE:0:1:2}, i.e., $\beta\in (0,2)$ and $\alpha=1$ and $\omega_{\nu}(\theta) \equiv 1$, has been also studied, see e.g., \cite[Section 7]{HiTa21} and \cite{HiSi24}. Nevertheless, $\alpha \in (0,1)$ is a case not previously investigated in the non-Gaussian framework for the Hardy equation, to the best of our knowledge.

It is worth mentioning that an extensive literature is devoted to the study of the heat equation with \textit{singular potentials}. For instance, one of the first papers on this topic was \cite{BaGo84} in the Euclidean setting. Recently, other studies have been developed in the setting of graded Lie groups; see, e.g. \cite{CRT22} where this theory is well documented. In this way, the present work provides some extensions to equations with memory, whose implementation is achieved by the fractional derivative in time of order $0<\alpha <1$.

The paper is organized as follows. In Section \ref{Sec:2} we review the standard facts on the fundamental solutions $Z,Y$ and we introduce some definitions. We also derive new properties of such solutions. The main result of this section is Lemma \ref{descomposition}. In the third section, we are concerned with the existence of local solutions and the main result is given in Theorem \ref{LocalResult}. Section \ref{Sec:4} establishes sufficient conditions for existence of global solutions. The main results are Theorems \ref{globalsolutionsmallu0} and \ref{Global2}, however, Corollary \ref{corol_global} provides a setting for existence of non-negative solutions. In section \ref{Sec:5}, an asymptotic analysis for solutions is given via the kernels $Z,Y$, which is stated in Theorem \ref{decay}. Section \ref{Sec:6} is devoted to the study of non-existence of positive solutions, the main result being given by Theorem \ref{Theorem_blow-up}. In the final section, we make a brief comparison of our results with the classical ones.
\section{Preliminaries}
\label{Sec:2}
In this section we compile certain properties of $Z$ and $Y$ from \cite{SoVe25, SoVe23} and we prove additional properties of them, which we later utilize throughout the paper. Since our methods work well with so-called \textit{mild solutions}, we introduce the non-linear operator
\begin{equation}
\label{OperatorM}
\mathcal{M}v(t) := S(t)u_0 + \int_{0}^{t} R(t-s)\left(|\cdot|^{-\gamma}|v|^{p-1}v(s)\right) ds
\end{equation}
on a set of continuous functions $v$ from $I\subset [0,\infty)$ to $L_q(\RR^d)$.
\begin{definition}
\label{DefSolution}
Let $1\leq q<\infty$ and $u_0\in L_q(\RR^d)$. The function $u\in C([0,T];L_q(\RR^d))$ is a {\bf local mild solution} of \eqref{FE:0} if it is measurable and satisfies $u=\mathcal{M}u$ for almost every $t\in [0, T]$. If $u\in C((0,\infty);L_q(\RR^d))$, we say that the solution is {\bf global}.
\end{definition}

\begin{proposition}\label{propo:Z:Y}
	Let $d\in \NN$, $\alpha\in (0,1)$, and $\beta\in (0,2)$. Under assumption $(H2)$, the following properties hold:
	\begin{itemize}
		\item[(i)] $\int_{\RR^d} Z(t,x)dx =1$ for all $t>0$, and $Z(t,x) = t^{-\frac{\alpha d}{\beta}} Z(1, t^{-\frac{\alpha}{\beta}} x)$ for all $t>0$ and $x\in \RR^d$.
		\item[(ii)] $\int_{\RR^d} Y(t,x)dx =g_{\alpha}(t)$ for all $t>0$, and $Y(t,x) = t^{-\frac{\alpha d}{\beta} + \alpha -1} Y(1, t^{-\frac{\alpha}{\beta}} x)$ for all $t>0$ and $x\in \RR^d$.
	\end{itemize}
\end{proposition}
For the proof, refer to \cite[Lemma 2.12]{SoVe22}. The first part of (ii) can be found in the proof of Theorem 2.14 in \cite{SoVe22}. Other properties of $Z$ and $Y$ are gathered from \cite[Theorem 2.8 and Theorem 2.10]{SoVe22}, respectively, in the following proposition.

\begin{proposition}\label{propo:lp:Z:Y}
		Let $d\in \NN$, $\alpha\in (0,1)$, and $\beta\in (0,2)$. Under assumption $(H2)$, the following properties are established:
		\begin{itemize}
			\item[(i)] The kernel $Z(t,\cdot)$ is in $L_r(\mathbb{R}^d)$ for all $t>0$ if and only if $1\leq r<\kappa_1$, where
			\[
			\kappa_1=\kappa_1(d,\beta):=
			\begin{cases}
				\frac{d}{d-\beta} & \text{if } d>\beta,\\
				\infty & \text{otherwise}.
			\end{cases}
			\]
			Moreover, the two-sided estimate 
			\begin{equation}
				\label{cotasZp}
				\norm{Z(t,\cdot)}_r\asymp t^{-\frac{\alpha d}{\beta}\left(1-\frac{1}{r}\right)},\; t>0
			\end{equation}
			holds for every $1\leq r<\kappa_1$. In the case of $d<\beta$, \eqref{cotasZp} remains true for $r=\infty$.\\
			\item[(ii)] The kernel $Y(t,\cdot)$ is in $L_r(\mathbb{R}^d)$ for all $t>0$ if and only if $1\leq r<\kappa_2$, where
			\[
			\kappa_2=\kappa_2(d,\beta):=
			\begin{cases}
				\frac{d}{d-2\beta} & \text{if } d>2\beta,\\
				\infty & \text{otherwise}.
			\end{cases}
			\]
			Moreover, the two-sided estimate 
			\begin{equation}
				\label{cotasYp}
				\norm{Y(t,\cdot)}_r\asymp t^{-\frac{\alpha d}{\beta}\left(1-\frac{1}{r}\right)+(\alpha-1)},\; t>0
			\end{equation}
			holds for every $1\leq r<\kappa_2$. In the case of $d<2\beta$, \eqref{cotasYp} remains true for $r=\infty$.
		\end{itemize}
\end{proposition}
The following property is proved in \cite[Proposition 2.3]{SoVe25}.
\begin{proposition}{\label{propo:z:phi}}
Let $d\in \NN$, $\alpha\in (0,1)$, and $\beta\in (0,2)$. Suppose that assumption $(H2)$ holds. Fix $\eta\in (0, d)$ and $R>1$. Let $u_0\in L_1(\RR^d)$ be the non-negative, radially symmetric function given by
	\[
	u_0(x) := |x|^{-\eta}1_{B(0,R)}(x) :=
	\begin{cases}
		|x|^{-\eta}, & |x|\leq R,\\
		0, & |x| >R.
	\end{cases}
	\]
	Let $z(t,x):=S(t)u_0(x)$ and $M:=\min\left\{z(t,\hat{x}) \, : \,  \hat{x}\in \mathbb{S}^{d-1}, \, 0\leq t\leq 1\right\}$. If $\rho\in (0,\alpha/\beta)$ then for any $\phi\geq M$ we have
	\[
	z(t,x) \geq \phi \, \text{ for } t^{\alpha/\beta} \leq |x| \leq t^{\rho}, \text{ for all } 0<t\leq \left(\frac{\phi}{M}\right)^{-\frac{1}{\eta \rho}}.
	\]
\end{proposition}

\begin{lemma}
Let $0<\gamma<d$. Under assumptions of Proposition \ref{propo:Z:Y}, the kernel $Y$ satisfies the estimate
\[
\int_{\RR^d}Y(t,x-y)|y|^{-\gamma}dy\leq Ct^{\alpha-1}|x|^{-\gamma}
\]
for all $t>0$ and $x\in\RR^d$, where the constant $C$ depends on $d,\alpha,\beta,\gamma$.
\end{lemma}
\begin{proof}
We define the sets $\Omega_1=\{y\in\RR^d: |y|\leq\frac{|x|}{2}\}$ and $\Omega_2=\{y\in\RR^d: |y|\geq\frac{|x|}{2}\}$. It follows that
\[
\int_{\RR^d}Y(t,x-y)|y|^{-\gamma}dy=\int_{\Omega_1}Y(t,x-y)|y|^{-\gamma}dy+\int_{\Omega_2}Y(t,x-y)|y|^{-\gamma}dy.
\]
Over $\Omega_2$, the integral is bounded by $|x|^{-\gamma}$. The fact that $\int_{\Omega_2}Y(t,x-y)|y|^{-\gamma}dy\leq\int_{\RR^d}Y(t,x-y)|y|^{-\gamma}dy$
and Proposition \ref{propo:Z:Y} yield
\[
\int_{\Omega_2}Y(t,x-y)|y|^{-\gamma}dy\leq 2^{\gamma}g_\alpha(t)|x|^{-\gamma}.
\]
We recall that $\Omega=|x-y|^{\beta} t^{-\alpha}$ and using the estimates \eqref{y:doble:est}, we split the integral over $\Omega_1$ into two parts, considering $\Omega\geq 1$ and $\Omega\leq 1$, respectively. We see that the first part yields the estimate
\begin{align*}
\int_{\Omega_1\cap\{\Omega\geq 1\}}Y(t,x-y)|y|^{-\gamma}dy&\leq Ct^{2\alpha-1}\int_{\Omega_1\cap\{\Omega\geq 1\}}|x-y|^{-\beta-d}|y|^{-\gamma}dy\\
&=Ct^{\alpha-1}\int_{\Omega_1\cap\{\Omega\geq 1\}}t^{\alpha}|x-y|^{-\beta}|x-y|^{-d}|y|^{-\gamma}dy\\
&=Ct^{\alpha-1}\int_{\Omega_1\cap\{\Omega\geq 1\}}\Omega^{-1}|x-y|^{-d}|y|^{-\gamma}dy\\
&\leq C2^dt^{\alpha-1}|x|^{-d}\int_{\Omega_1}|y|^{-\gamma}dy
\end{align*}
since $|x-y|\geq |x|-|y|\geq |x|-\frac{|x|}{2}=\frac{|x|}{2}$. Therefore,
\[
\int_{\Omega_1\cap\{\Omega\geq 1\}}Y(t,x-y)|y|^{-\gamma}dy\lesssim t^{\alpha-1}|x|^{-\gamma}. 
\]
The last part is analysed according to the cases $d<2\beta$, $d=2\beta$ and $d>2\beta$, respectively, obtaining the estimate
\[
\int_{\Omega_1\cap\{\Omega\leq 1\}}Y(t,x-y)|y|^{-\gamma}dy\lesssim t^{\alpha-1}|x|^{-\gamma}. 
\]
\end{proof}

\begin{lemma}
\label{descomposition}
Let $d\in \NN$, $\alpha\in (0,1)$, and $\beta\in (0,2)$. Suppose that assumption $(H2)$ holds. Let $0<\gamma<d$, define $m:=\frac{d}{\gamma}$ and fix $\epsilon >0$ satisfying $m-\epsilon>1$. Suppose that there exists $1< q_1<\infty$ such that $r_1:=\frac{q_1(m-\epsilon)}{q_1+m-\epsilon}\geq 1$, $r_2:=\frac{q_1(m+\epsilon)}{q_1+m+\epsilon}\geq 1$. Define 
	\[
	\kappa_3 =\kappa_3(d,\beta,r_1)=
	\begin{cases}
		\frac{r_1d}{d - 2r_1\beta} & \text{if } d > 2r_1\beta, \\
		\infty & \text{otherwise}.
	\end{cases}
	\]
Let $q_1\leq q_2<\kappa_3$. Then the operator $R$ satisfies the estimate
\begin{align*}
\left|R(t)|\cdot|^{-\gamma}v\right|_{q_2}\leq Ct^{\alpha-1-\frac{\alpha d}{\beta}\left(\frac{1}{q_1}-\frac{1}{q_2}\right)-\frac{\alpha \gamma}{\beta}}|v|_{q_1}
\end{align*}
for all $t>0$ and $v\in L_{q_1}(\RR^d)$, where the constant $C$ depends on $d,\alpha,\beta,\gamma, \epsilon, q_1, q_2$.
\end{lemma}

\begin{proof}
Firstly, we note that $q_1>r_2>r_1$. Definition of $r_1,r_2$ implies that
\begin{equation}\label{Holder:exponents:1:2}
		\frac{1}{r_1} = \frac{1}{m - \epsilon} + \frac{1}{q_1} \quad \text{and} \quad \frac{1}{r_2} = \frac{1}{m + \epsilon} + \frac{1}{q_1}.
	\end{equation}
Therefore we can find numbers $l_1,l_2$ such that
\begin{equation}\label{Young:exponents:1:2}
		1 + \frac{1}{q_2} = \frac{1}{l_1} + \frac{1}{r_1} = \frac{1}{l_2} + \frac{1}{r_2}.
	\end{equation}
Conditions on $q_2$ guarantee that $1< l_1,l_2<\kappa_2$, where $\kappa_2$ is as in Proposition \ref{propo:lp:Z:Y} (ii). 

We split the function $|\cdot|^{-\gamma}$ as follows:
	\[
	|\cdot|^{-\gamma} = \varphi_1 + \varphi_2,
	\]
	where $\varphi_1 \in L_{m - \epsilon}(\mathbb{R}^d)$ and $\varphi_2 \in L_{m + \epsilon}(\mathbb{R}^d)$. Thus, we apply Young's inequality for the exponents \eqref{Young:exponents:1:2} and Hölder's inequalities for \eqref{Holder:exponents:1:2}, together with the estimate of $Y$ in \eqref{cotasYp}, to obtain
\begin{align*}
		\left| R(t) | \cdot |^{-\gamma} v \right|_{q_2} &\leq \left| Y(t, \cdot) \star (\varphi_1 v) \right|_{q_2} + \left| Y(t, \cdot) \star (\varphi_2 v) \right|_{q_2} \\
		&\leq |Y(t, \cdot)|_{l_1} |\varphi_1 v|_{r_1} + |Y(t, \cdot)|_{l_2} |\varphi_2 v|_{r_2} \\
		&\leq C_1 t^{-\frac{\alpha d}{\beta} \left( \frac{1}{r_1} - \frac{1}{q_2} \right) + (\alpha - 1)} |\varphi_1|_{m - \epsilon} |v|_{q_1} \\
		&\quad + C_2 t^{-\frac{\alpha d}{\beta} \left( \frac{1}{r_2} - \frac{1}{q_2} \right) + (\alpha - 1)} |\varphi_2|_{m + \epsilon} |v|_{q_1}.
	\end{align*}
	In particular, for $t = 1$, we have
	\begin{equation}
		\label{estimate1}
		\left| R(1) | \cdot |^{-\gamma} v \right|_{q_2} \leq C |v|_{q_1}.
	\end{equation}
	
Next, let $\lambda >0$. We define the dilation operator $\mathcal{D}_\lambda$ given by $(\mathcal{D}_\lambda\phi)(x)=\phi(\lambda x)$ for any $\phi\in L_q(\RR^d)$, $1\leq q<\infty$. We recall that $|\mathcal{D}_\lambda\phi|_q=\lambda^{-\frac{d}{q}}|\phi|_q$. This operator and Proposition \ref{propo:Z:Y} show that
\begin{align*}
Y(t,\cdot)\star(|\cdot|^{-\gamma} v)(x)&=\int_{\RR^d}Y(t,x-y)|y|^{-\gamma}v(y)dy\\
&=t^{-\frac{\alpha d}{\beta} + \alpha -1}\int_{\RR^d} Y(1, t^{-\frac{\alpha}{\beta}} (x-y))|y|^{-\gamma}v(y)dy\\
&=t^{\alpha -1}\int_{\RR^d} Y(1, t^{-\frac{\alpha}{\beta}}x-w)t^{-\frac{\alpha\gamma}{\beta}}|w|^{-\gamma}v(t^{\frac{\alpha}{\beta}}w)dw\\
&=t^{\alpha -1-\frac{\alpha\gamma}{\beta}}\int_{\RR^d} Y(1, t^{-\frac{\alpha}{\beta}}x-w)|w|^{-\gamma}(\mathcal{D}_{t^{\frac{\alpha}{\beta}}}v)(w)dw\\
&=t^{\alpha -1-\frac{\alpha\gamma}{\beta}}Y(1,\cdot)\star(|\cdot|^{-\gamma}\mathcal{D}_{t^{\frac{\alpha}{\beta}}}v)(t^{-\frac{\alpha}{\beta}}x).
\end{align*}
This implies that
\[
Y(t,\cdot)\star(|\cdot|^{-\gamma} v)=t^{\alpha -1-\frac{\alpha\gamma}{\beta}}D_{t^{-\frac{\alpha}{\beta}}}\left(Y(1,\cdot)\star(|\cdot|^{-\gamma}\mathcal{D}_{t^{\frac{\alpha}{\beta}}}v)\right).
\]
The proof is completed by referring to \eqref{estimate1}, because
\begin{align*}
\left|R(t)|\cdot|^{-\gamma}v\right|_{q_2}&=\left|Y(t,\cdot)\star(|\cdot|^{-\gamma} v)\right|_{q_2}\\
&=t^{\alpha -1-\frac{\alpha\gamma}{\beta}}\left|D_{t^{-\frac{\alpha}{\beta}}}\left(Y(1,\cdot)\star(|\cdot|^{-\gamma}\mathcal{D}_{t^{\frac{\alpha}{\beta}}}v)\right)\right|_{q_2}\\
&=t^{\alpha -1-\frac{\alpha\gamma}{\beta}+\frac{\alpha d}{\beta q_2}}\left|Y(1,\cdot)\star(|\cdot|^{-\gamma}\mathcal{D}_{t^{\frac{\alpha}{\beta}}}v)\right|_{q_2}\\
&\leq C t^{\alpha -1-\frac{\alpha\gamma}{\beta}+\frac{\alpha d}{\beta q_2}}\left|\mathcal{D}_{t^{\frac{\alpha}{\beta}}}v)\right|_{q_1}\\
&=C t^{\alpha -1-\frac{\alpha\gamma}{\beta}+\frac{\alpha d}{\beta q_2}-\frac{\alpha d}{\beta q_1}}|v|_{q_1}.
\end{align*}
\end{proof}
Before finishing this section, we define an Osgood-type function $f$ whose properties we will use later. For more details on this type of function, see \cite[Section 3]{LaRoSi13}. 

Considering the given parameter $p>1$, we fix a number $\phi_0>2^{\frac{1}{p-1}}$. Define a sequence $(\phi_i)_{i\in \NN}$ as
\begin{align}
\label{phi_i}
\phi_{i}=\phi_{i-1}^p,\quad i\geq 1,
\end{align}
such that
$$1<\phi_{i-1}<\frac{\phi_i}{2},\quad i\geq 1, $$
and 
$\phi_i\rightarrow\infty$ as $i\rightarrow\infty$. The function $f$ is defined as
\begin{equation}
\label{fOsgood}
f(s) :=\begin{cases}\left(1-\phi_0^{1-p}\right)s^p & \text{if }s\in J_0=[0,\phi_0],\\ \phi_i-\phi_{i-1} & \text{if }s\in I_i=\left[\phi_{i-1},\frac{\phi_i}{2}\right],~i\geq 1,\\l_i(s) &\text{if }s\in J_i=\left(\frac{\phi_i}{2},\phi_i\right),~i\geq 1,  \end{cases}
\end{equation}
where $l_i$ is a polynomial of order $1$ that interpolates between the values of $f$ at $\frac{\phi_i}{2}$ and $\phi_i$. It is not difficult to check that $f\geq \frac{\phi_{i+1}}{2}$ on $[\phi_i,\infty)$, $i\geq 1$.

\section{Local well-posedness}
\label{Sec:3}
In this section we show that the Hardy non-Gaussian equation \eqref{FE:0} is locally solvable in the Banach space $Y_T=C([0,T];L_q(\RR^d))$ with the usual norm $\norm{v}_{Y_T}=\sup_{0\leq t\leq T}|v(t)|_q$ and $p<q<\infty$, according to Definition \ref{DefSolution}. That is, we show that the operator
\[
\mathcal{M}v(t)= S(t)u_0 + \int_{0}^{t} R(t-s)\left(|\cdot|^{-\gamma}|v|^{p-1}v(s)\right) ds
\]
given by \eqref{OperatorM}, has a unique fixed point in $Y_T$ for some $T>0$, where  $u_0\in L_q(\RR^d)$ and the operators $S$, $R$ are given in \eqref{operatosRS}. From Proposition \ref{propo:Z:Y} it is straightforward to see that
\[
|S(t)u_0|_q=|Z(t,\cdot)\star u_0|_q\leq |u_0|_q.
\]
For the non-linear term of $\mathcal{M}$ we apply Lemma \ref{descomposition}, considering $q_2=q$ and $q_1=\frac{q}{p}$. It follows that
\begin{align*}
\left|\int_{0}^{t} R(t-s)|\cdot|^{-\gamma}|v|^{p-1}v(s)ds\right|_q&\leq \int_{0}^{t} |R(t-s)|\cdot|^{-\gamma}|v|^{p-1}v(s)|_q ds\\
&\leq C\int_{0}^{t}(t-s)^{\alpha-1-\frac{\alpha d}{\beta}\left(\frac{p}{q}-\frac{1}{q}\right)-\frac{\alpha \gamma}{\beta}}||v|^{p-1}v(s)|_{\frac{q}{p}}ds
\end{align*}
and the fact that $||v|^{p-1}v(s)|_{\frac{q}{p}}=|v(s)|_q^{p}$ imply
\begin{align*}
\left|\int_{0}^{t} R(t-s)|\cdot|^{-\gamma}|v|^{p-1}v(s)ds\right|_q&\leq C\norm{v}_{Y_T}^{p}\int_{0}^{t}(t-s)^{\alpha-1-\frac{\alpha d}{\beta}\left(\frac{p}{q}-\frac{1}{q}\right)-\frac{\alpha \gamma}{\beta}}ds\\
&=C_1\norm{v}_{Y_T}^{p}t^{\alpha-\frac{\alpha d}{\beta q}\left(p-1\right)-\frac{\alpha \gamma}{\beta}}
\end{align*}
whenever $q>\frac{d}{\beta-\gamma}(p-1)$ and $\gamma<\beta$. Hence, it makes sense to define
\begin{equation}
\label{Cexponent}
q_c:=\frac{d}{\beta-\gamma}(p-1)
\end{equation}
as the \textbf{critical exponent} for the local well-posedness. The continuity of $t\mapsto \mathcal{M}v$ is readily checked using similar arguments as in \cite[Section 3]{SoVe22}. This shows that $\mathcal{M}$ is well defined.
\begin{remark}
Critical exponents are usually defined grater than $1$, however, in this section we only need that $q>q_c$ with $q_c$ given by \eqref{Cexponent}.
\end{remark}
Next, we shall prove that $\mathcal{M}$ is a contraction on the closed ball $B=\{v\in Y_T: |v|_{Y_T}\leq 2 \norm{u_0}_{q}\}$. Indeed, let $v,w\in B$. It follows that
\begin{align*}
|\mathcal{M}v(t)|_q&\leq |S(t)u_0|_q + \int_{0}^{t} |R(t-s)|\cdot|^{-\gamma}|v|^{p-1}v(s)|_q ds\\
&\leq |u_0|_q+C\int_{0}^{t}(t-s)^{\alpha-1-\frac{\alpha d}{\beta}\left(\frac{p}{q}-\frac{1}{q}\right)-\frac{\alpha \gamma}{\beta}}||v|^{p-1}v(s)|_{\frac{q}{p}}ds\\
&\leq |u_0|_q+C_1\norm{v}_{Y_T}^{p}T^{\alpha-\frac{\alpha d}{\beta q}\left(p-1\right)-\frac{\alpha \gamma}{\beta}}\\
&\leq |u_0|_q+C_12^p|u_0|_q^{p}T^{\alpha-\frac{\alpha d}{\beta q}\left(p-1\right)-\frac{\alpha \gamma}{\beta}}
\end{align*}
which implies that $|\mathcal{M}v|_{Y_T}\leq 2|u_0|_q$ for $T$ sufficiently small. Besides,
\begin{align*}
|\mathcal{M}v(t)-\mathcal{M}w(t)|_q&\leq \int_{0}^{t} |R(t-s)|\cdot|^{-\gamma}\left(|v|^{p-1}v(s)-|w|^{p-1}w(s)\right)|_q ds\\
&\leq C\int_{0}^{t}(t-s)^{\alpha-1-\frac{\alpha d}{\beta}\left(\frac{p}{q}-\frac{1}{q}\right)-\frac{\alpha \gamma}{\beta}}||v|^{p-1}v(s)-|w|^{p-1}w(s)|_{\frac{q}{p}}ds\\
&\leq C|v+w|^{p-1}_{Y_T}|v-w|_{Y_T}\int_{0}^{t}(t-s)^{\alpha-1-\frac{\alpha d}{\beta q}\left(p-1\right)-\frac{\alpha \gamma}{\beta}}ds
\end{align*}
where the last estimate is obtained by applying the property
\begin{equation}
\label{PropiedadResta}
|\,|a|^c a-|b|^c b\,|\lesssim |\,a-b\,|(|a|^c+|b|^c)\lesssim |\,a-b\,|(\,|a|+|b|\,)^c,\quad a, b\in\RR,~c>0,
\end{equation}
and Hölder's inequality. We see that  
\[
|\mathcal{M}v(t)-\mathcal{M}w(t)|_q\leq C_1 (4|u_0|_q)^{p-1}|v-w|_{Y_T}T^{\alpha-\frac{\alpha d}{\beta q}\left(p-1\right)-\frac{\alpha \gamma}{\beta}}
\]
and therefore $\mathcal{M}$ is a contraction for small enough $T$. From the contraction mapping principle, it follows that $\mathcal{M}$ has a unique fixed point $u$ in $B$. Finally, we show that the solution $u$ is unique in $Y_T$. For this purpose we suppose that there exists another function $\tilde{u}\in C([0,T];L_q(\RR^d))$ satisfying  $\mathcal{M}\tilde{u} =\tilde{u}$.
Let $t\in (0,T]$. Previous work shows that
\begin{align*}
\norm{u(t)-\tilde{u}(t)}_q &= \norm{\mathcal{M}u(t)-\mathcal{M}\tilde{u}(t)}_q\\
&\leq C|u+\tilde{u}|^{p-1}_{Y_T}\int_{0}^{t}(t-s)^{\alpha-\frac{\alpha d}{\beta q}\left(p-1\right)-\frac{\alpha \gamma}{\beta}-1}|u(s)-\tilde{u}(s)|_{q}ds.
\end{align*}
Since the condition $q>q_c$ ensures that $\alpha-\frac{\alpha d}{\beta q}\left(p-1\right)-\frac{\alpha \gamma}{\beta}>0$, the uniqueness in $Y_T$ follows from \cite[Theorem 1]{YGD07}. In this way, we have proved the following result.
\begin{theorem}
\label{LocalResult}
Let $0<\gamma<\min(\beta,d)$ and $\max(p,q_c)<q<\infty$. Assume that Lemma \ref{descomposition} is satisfied with $q_1=\frac{q}{p}$ and $q_2=q$. Let $u_0\in L_{q}(\mathbb{R}^d)$. If $T$ is sufficiently small then the Cauchy problem \eqref{FE:0} is well-posedness and it has a unique solution $u\in C([0,T];L_q(\RR^d))$.
\end{theorem}

\section{Global Solutions}
\label{Sec:4}
In this section we consider the same critical exponent $q_c$ given by \eqref{Cexponent}, but in such a way that the solution is global. 
\begin{theorem}
	\label{globalsolutionsmallu0}
Let $0<\gamma<\min(\beta,d)$ and $q_c\geq \max\left(1, \frac{d}{\beta}\right)$. Suppose that $\alpha p>p-1$ and that $\max\left(p,q_c\right)<q<\frac{\alpha dp(p-1)}{\beta(\alpha p -p+1)}$. Assume that Lemma \ref{descomposition} is satisfied with $q_1=\frac{q}{p}$ and $q_2=q$. Let $u_0\in L_{q_c}(\mathbb{R}^d)$. If $|u_0|_{q_c}$ is sufficiently small, then the problem \eqref{FE:0} has a global solution $u\in C((0,\infty);L_q(\mathbb{R}^d))$. Moreover, $|u(t)|_q\rightarrow 0$ as $t\rightarrow \infty$.
\end{theorem}
\begin{proof}
We define the Banach space
\[
E:=C((0,\infty);L_q(\RR^d))
\]
with the norm
	\[
	\left\|v\right\|_{E}:=\esssup_{t> 0}t^{\frac{\alpha d}{\beta}\left(\frac{1}{q_c}-\frac{1}{q}\right)}\norm{v(t,\cdot)}_{q}.
	\]
On this space, again, we consider the operator $\mathcal{M}$ given by \eqref{OperatorM},
\[
\mathcal{M}v(t)= S(t)u_0 + \int_{0}^{t} R(t-s)\left(|\cdot|^{-\gamma}|v|^{p-1}v(s)\right) ds.
\]
Using properties of the kernel $Z$ from Proposition \ref{propo:lp:Z:Y}, Lemma \ref{descomposition} considering $q_2=q$ and $q_1=\frac{q}{p}$, together with Young's convolution inequality applied to the relationships $1+ \frac{1}{q} = \frac{1}{q_c} + \frac{1}{p}$ and $1+ \frac{1}{q} = \frac{k}{q} + \frac{1}{p}$ respectively, it follows that this operator is well-defined. Indeed, for $v\in E$ we have that
\begin{align*}
\norm{\mathcal{M}v(t)}_q & \leq \norm{S(t)u_0}_q + \int_{0}^{t} |R(t-s)|\cdot|^{-\gamma}|v|^{p-1}v(s)|_q ds\\
&= \norm{Z(t)\star u_0}_q + \int_{0}^{t} |R(t-s)|\cdot|^{-\gamma}|v|^{p-1}v(s)|_q ds\\
&\leq C_1 t^{-\frac{\alpha d}{\beta}\left(\frac{1}{q_c}-\frac{1}{q}\right)}\norm{u_0}_{q_c}+
C\int_{0}^{t}(t-s)^{\alpha-1-\frac{\alpha d}{\beta}\left(\frac{p}{q}-\frac{1}{q}\right)-\frac{\alpha \gamma}{\beta}}||v|^{p-1}v(s)|_{\frac{q}{p}}ds.
\end{align*}
As in previous section
\[
\norm{|v|^{p-1}v(s)}_{\frac{q}{p}}=\norm{v(s)}_q^p
\]
and thus
\begin{align*}
\norm{\mathcal{M}v(t)}_q&\leq C_1 t^{-\frac{\alpha d}{\beta}\left(\frac{1}{q_c}-\frac{1}{q}\right)}\norm{u_0}_{q_c}+ C\left\|v\right\|_{E}^{p}\int_{0}^{t}(t-s)^{\alpha-1-\frac{\alpha d}{\beta}\left(\frac{p}{q}-\frac{1}{q}\right)-\frac{\alpha \gamma}{\beta}}s^{-\frac{\alpha d p}{\beta}\left(\frac{1}{q_c}-\frac{1}{q}\right)}ds\\
&=C_1t^{-\frac{\alpha d}{\beta}\left(\frac{1}{q_c}-\frac{1}{q}\right)}\norm{u_0}_{q_c}\\
&~~~~ + C\left\|v\right\|_E^p t^{\alpha-\frac{\alpha d (p-1)}{\beta q}-\frac{\alpha \gamma}{\beta}-\frac{\alpha d p}{\beta}\left(\frac{1}{q_c}-\frac{1}{q}\right)}\int_{0}^{1} (1-\tau)^{\alpha-\frac{\alpha d (p-1)}{\beta q}-\frac{\alpha \gamma}{\beta}-1}\tau^{-\frac{\alpha d p}{\beta}\left(\frac{1}{q_c}-\frac{1}{q}\right)}d\tau.
\end{align*}
Definition of $q_c$ and conditions on $q$, show that
\begin{align*}
\norm{\mathcal{M}v(t)}_q&\leq C_1 t^{-\frac{\alpha d}{\beta}\left(\frac{1}{q_c}-\frac{1}{q}\right)}\norm{u_0}_{q_c}+C\left\|v\right\|_E^p t^{\alpha+\frac{\alpha d}{\beta q}-\frac{\alpha d p}{\beta q_c}-\frac{\alpha \gamma}{\beta}}\frac{\Gamma\left(\alpha-\frac{\alpha d(p-1)}{\beta q}-\frac{\alpha \gamma}{\beta}\right)\Gamma\left(1-\frac{\alpha d p}{\beta}\left(\frac{1}{q_c}-\frac{1}{q}\right)\right)}{\Gamma\left(1+\alpha+\frac{\alpha d}{\beta q}-\frac{\alpha d p}{\beta q_c}-\frac{\alpha \gamma}{\beta}\right)}\\
&=C_1t^{-\frac{\alpha d}{\beta}\left(\frac{1}{q_c}-\frac{1}{q}\right)}\norm{u_0}_{q_c}+C\left\|v\right\|_E^p t^{-\frac{\alpha d}{\beta}\left(\frac{1}{q_c}-\frac{1}{q}\right)}\frac{\Gamma\left(\alpha-\frac{\alpha d(p-1)}{\beta q}-\frac{\alpha \gamma}{\beta}\right)\Gamma\left(1-\frac{\alpha d p}{\beta}\left(\frac{1}{q_c}-\frac{1}{q}\right)\right)}{\Gamma\left(1-\frac{\alpha d}{\beta}\left(\frac{1}{q_c}-\frac{1}{q}\right)\right)}.
\end{align*}
This proves that $\norm{\mathcal{M}v(t)}_q$ exists for all $t>0$. Furthermore, by defining $C_2:=C\frac{\Gamma\left(\alpha-\frac{\alpha d(p-1)}{\beta q}-\frac{\alpha \gamma}{\beta}\right)\Gamma\left(1-\frac{\alpha d p}{\beta}\left(\frac{1}{q_c}-\frac{1}{q}\right)\right)}{\Gamma\left(1-\frac{\alpha d}{\beta}\left(\frac{1}{q_c}-\frac{1}{q}\right)\right)}$, we obtain
\begin{align*}
t^{\frac{\alpha d}{\beta}\left(\frac{1}{q_c}-\frac{1}{q}\right)}\norm{\mathcal{M}v(t)}_q&\leq C_1\norm{u_0}_{q_c}+ C_2\left\|v\right\|_E^p,
\end{align*}
which consequently establishes $\left\|\mathcal{M}v\right\|_E<\infty$. The continuity of $\mathcal{M}v$ on $(0,\infty)$, as in previous section, follows from \cite[Section 3]{SoVe22}. This shows that $\mathcal{M}$ is well defined.

Next, we consider the closed ball $B=\{v\in E: \left\|v\right\|_E\leq 2 C_1\norm{u_0}_{q_c}\}$. Again, we shall demonstrate that $\mathcal{M}$ has a unique fixed point in $B$ by utilizing the contraction mapping principle. The preceding work establishes that if $v\in B$ then
\begin{align*}
\left\|\mathcal{M}v\right\|_E\leq C_1\norm{u_0}_{q_c}+C_2 2^p C_1^p\norm{u_0}_{q_c}^p \leq \left( C_1+ C_2 2^pC_1^p \norm{u_0}_{q_c}^{p-1}\right)\norm{u_0}_{q_c}
\end{align*}
and thus $\mathcal{M}v \in B$ for small enough $\norm{u_0}_{q_c}$. For $v,w\in B$, using the property \eqref{PropiedadResta} and Hölder's inequality, we see that
\begin{align*}
\norm{\mathcal{M}v(t)-\mathcal{M}w(t)}_{q}&\leq C\left\|v+w\right\|_E^{p-1}\left\|v-w\right\|_E\int_{0}^{t} (t-s)^{\alpha-\frac{\alpha d(p-1)}{\beta q}-\frac{\alpha\gamma}{\beta}-1}s^{-\frac{\alpha d p}{\beta}\left(\frac{1}{q_c}-\frac{1}{q}\right)}ds\\
&=C_2\left\|v+w\right\|_E^{p-1}\left\|v-w\right\|_E t^{-\frac{\alpha d}{\beta}\left(\frac{1}{q_c}-\frac{1}{q}\right)}
\end{align*}
which implies that
\[
\left\|\mathcal{M}v-\mathcal{M}w\right\|_E\leq C_3\norm{u_0}_{q_c}^{p-1}\left\|v-w\right\|_E.
\]
As a result, $\mathcal{M}$ is a contraction whenever $\norm{u_0}_{q_c}$ is sufficiently small. Consequently, there exists a unique fixed point $u\in B$ of $\mathcal{M}$, that is, $\mathcal{M}u = u$ and $\norm{u(t,\cdot)}_{q}\leq 2C_1|u_0|_qt^{-\frac{\alpha d}{\beta}\left(\frac{1}{q_c}-\frac{1}{q}\right)}$ for all $t>0$. 
\end{proof}
\begin{remark}
Unlike the local case, conditions of Theorem \ref{globalsolutionsmallu0} does not allow a proof of the uniqueness in the space $E$ due to the presence of a non-integrable singularity in $s$.
\end{remark}
\begin{theorem}
\label{Global2}
Let $0<\gamma<\min(\beta,d)$ and $q_c > \max\left(1, \frac{d}{\beta}\right)$. Suppose that $\alpha p>p-1$ and that $\max\left(p,q_c\right)<q<\frac{\alpha dp(p-1)}{\beta(\alpha p -p+1)}$. Assume that Lemma \ref{descomposition} is satisfied with $q_1=\frac{q}{p}$ and $q_2=q$. Let $u_0\in L_{1,loc}(\RR^d)$ such that $|u_0|\leq A |\cdot|^{-\frac{d}{q_c}}$. If the constant $A$ is sufficiently small, then the problem \eqref{FE:0} has a global solution $u\in C((0,\infty);L_q(\mathbb{R}^d))$ with $|u(t)|_q\rightarrow 0$ as $t\rightarrow \infty$.
\end{theorem}
\begin{proof}
We denote by $B(\epsilon)$ the closed ball in $\RR^d$, with center at the origin and radius $\epsilon >0$. As in \cite[Section 4]{STW17} we write 
\[
|\cdot|^{-\frac{d}{q_c}}=|\cdot|^{-\frac{\beta-\gamma}{p-1}}=\phi_1+\phi_2,
\]
where $\phi_1(x)=\begin{cases}|x|^{-\frac{\beta-\gamma}{p-1}} &\text{if } x\in B(1),\\~0 &\text{if } x\in \RR^d\setminus B(1),  \end{cases}$ and $\phi_2(x)=\begin{cases}~0 &\text{if } x\in B(1),\\|x|^{-\frac{\beta-\gamma}{p-1}} &\text{if } x\in \RR^d\setminus B(1).  \end{cases}$

Definition of $q_c$, given in \eqref{Cexponent}, and the hypothesis $\max\left(1, \frac{d}{\beta}\right)<q_c<q$, guarantee that $\phi_1\in L_{s_1}(\RR^d)$ for any $\frac{d}{\beta}< s_1 <q_c$ and $\phi_2\in L_{s_2}(\RR^d)$ for any $q_c < s_2 < q $. It follows that 
\begin{align*}
|Z(t,\cdot)\star|\cdot|^{-\frac{\beta-\gamma}{p-1}}|_q&= |Z(t)\star(\phi_1+\phi_2)|_q\\
& \leq |Z(t)\star\phi_1|_{q}+|Z(t)\star\phi_2|_{q}\\
& \leq C_1t^{-\frac{\alpha d}{\beta}\left(\frac{1}{s_1}-\frac{1}{q}\right)}|\phi_1|_{s_1}+C_2t^{-\frac{\alpha d}{\beta}\left(\frac{1}{s_2}-\frac{1}{q}\right)}|\phi_2|_{s_2}
\end{align*}
and therefore
\[
|Z(1,\cdot)\star|\cdot|^{-\frac{\beta-\gamma}{p-1}}|_q\leq C,
\]
where $C$ is a  constant that depends on $d,\beta,\alpha,\gamma,p, s_1, s_2$. Now, as in proof of Lemma \ref{descomposition}, we derive that
\[
|Z(t,\cdot)\star u_0|\leq A t^{-\frac{\alpha(\beta-\gamma)}{\beta(p-1)}}D_{t^{-\frac{\alpha}{\beta}}}\left(Z(1,\cdot)\star|\cdot|^{-\frac{\beta-\gamma}{p-1}}\right).
\]
This implies that
\begin{align*}
|S(t)u_0|_q&\leq A t^{-\frac{\alpha(\beta-\gamma)}{\beta(p-1)}}\left|D_{t^{-\frac{\alpha}{\beta}}}\left(Z(1,\cdot)\star|\cdot|^{-\frac{\beta-\gamma}{p-1}}\right)\right|_q\\
&=A t^{-\frac{\alpha(\beta-\gamma)}{\beta(p-1)}}t^{\frac{\alpha d}{\beta q}}\left|Z(1,\cdot)\star|\cdot|^{-\frac{\beta-\gamma}{p-1}}\right|_q\\
&\leq A C t^{-\frac{\alpha(\beta-\gamma)}{\beta(p-1)}}t^{\frac{\alpha d}{\beta q}}\\
&= A C t^{-\frac{\alpha d}{\beta q_c}}t^{\frac{\alpha d}{\beta q}}\\
&= A C t^{-\frac{\alpha d}{\beta}\left(\frac{1}{q_c}-\frac{1}{q}\right)}.
\end{align*}
Finally, the existence of a solution $u\in C((0,\infty);L_q(\RR^d))$ follows from the contraction mapping principle, as in the proof of Theorem \ref{globalsolutionsmallu0}, but considering the closed ball $B=\{v\in E: \left\|v\right\|_E\leq 2 AC\}$.
\end{proof}
We conclude this section with a simple result for non-negative solutions.
\begin{corollary}
\label{corol_global}
Under assumptions of Theorem \ref{Global2}, there exists a non-negative global solution whenever $u_0\geq 0$.
\end{corollary}
\begin{proof}
The proof of Theorem \ref{Global2} ensures the existence of non-negative fixed points, by considering the closed set $B=\{v\in E: \left\|v\right\|_E\leq 2 AC\text{ and } v\geq 0\}$.
\end{proof}

\section{Asymptotic behaviour of solutions}
\label{Sec:5}
We prove an optimal behaviour which describes how the solutions, belonging to the Banach space $E$ defined in previous section, decay in time in terms of the fundamental solutions to equation \eqref{FE:0}. For this purpose we need the following assumption on the spectral measure:
\begin{itemize}
	\item[(H3)] The spectral measure $ \nu$ has a strictly positive density, such that the function $\omega_{\nu}$ in \eqref{omega:nu} is strictly positive and differentiable $(d + 2 + [\beta])$ times on $\mathbb{S}^{d-1}$.
\end{itemize}
Again, critical exponent $q_c$ defined by \eqref{Cexponent} plays an important role.
\begin{theorem}
\label{decay}
Let $\alpha\in (0,1)$, $\beta\in (1,2)$, $p>1$ and $0<\gamma<\min(\beta,d)$. Assume the hypothesis $(H3)$ holds. Suppose that $q_c\geq 1$, that $\alpha p>p-1$, that $\max\left(p,q_c\right)<q<\min\left(\frac{\alpha dp(p-1)}{\beta(\alpha p -p+1)},\kappa_1\right)$ and that $\frac{q}{q-p}<\frac{d}{\gamma}$, with $\kappa_1$ as in Proposition \ref{propo:lp:Z:Y}. Suppose that there exist numbers $r,l$ satisfying
\begin{equation*}
\label{pqr}
1\leq r<\frac{d}{d-1},\quad 1\leq l<\frac{d}{d+1-\beta},\quad\frac{1}{q}+1=\frac{1}{l}+\frac{1}{r},
\end{equation*}
such that the initial data $u_0\in L_1(\RR^d)$ and $\norm{\cdot}u_0\in L_r(\RR^d)$. If $u\in C((0,\infty);L_q(\RR^d))$ is a global solution to the Cauchy problem \eqref{FE:0}, then $u$ has the asymptotic behaviour
\[
\norm{u(t,\cdot)-AZ(t,\cdot)-BY(t,\cdot)}_q\rightarrow 0
\]
as $t\rightarrow\infty$, with the constants
\[
A=\displaystyle\int_{\RR^d} u_0(y)dy
\]
and
\[
B=\displaystyle\int_0^\infty\int_{\RR^d}|y|^{-\gamma}|u(s,y)|^{p-1}u(s,y)dyds.
\]
\end{theorem}
\begin{proof}

\begin{align*}
\norm{u(t,\cdot)-AZ(t,\cdot)-BY(t,\cdot)}_q&=\norm{S(t)u_0 + \int_{0}^{t} R(t-s)\left(|\cdot|^{-\gamma}|u|^{p-1}u(s)\right) ds-AZ(t,\cdot)-BY(t,\cdot)}_q\\
&\leq\norm{S(t)u_0-AZ(t,\cdot)}_q+\norm{\int_{0}^{t} R(t-s)\left(|\cdot|^{-\gamma}|u|^{p-1}u(s)\right) ds-BY(t,\cdot)}_q.
\end{align*}
The estimate for $\norm{S(t)u_0-AZ(t,\cdot)}_q$ follows straightforwardly from \cite[Lemmata 6.1-6.2]{SoVe22}, producing
\begin{equation}
\label{decayinit}
\norm{S(t)u_0-AZ(t,\cdot)}_q\lesssim \norm{\nabla Z(t,\cdot)}_{l}\norm{\norm{\cdot}u_0}_r\lesssim t^{-\frac{\alpha d}{\beta}\left(\frac{1}{r}-\frac{1}{q}\right)-\frac{\alpha}{\beta}}\norm{\norm{\cdot}u_0}_r.
\end{equation}
For the non linear estimate, we obtain that
\begin{align*}
&\norm{\int_{0}^{t} R(t-s)\left(|\cdot|^{-\gamma}|u|^{p-1}u(s)\right) ds-BY(t,\cdot)}_q\\
&\leq\norm{\int_{0}^{\frac{t}{2}}R(t-s)\left(|\cdot|^{-\gamma}|u|^{p-1}u(s)\right) ds}_q+\norm{\int_{\frac{t}{2}}^tR(t-s)\left(|\cdot|^{-\gamma}|u|^{p-1}u(s)\right) ds}_q+|B|\norm{Y(t,\cdot)}_q\\
&\leq\int_{0}^{\frac{t}{2}}\norm{Y(t-s,\cdot)\star\left(|\cdot|^{-\gamma}|u|^{p-1}u(s)\right)}_q ds+\int_{\frac{t}{2}}^t\norm{Y(t-s,\cdot)\star\left(|\cdot|^{-\gamma}|u|^{p-1}u(s)\right)}_q ds+|B|\norm{Y(t,\cdot)}_q\\
&=:J_1(t)+J_2(t)+J_3(t).
\end{align*}
To estimate both $J_1$ and $J_2$, we use the hypotheses $\frac{q}{q-p}<\frac{d}{\gamma}$ and $p<q<\kappa_1$, along with Hölder's, Young's and Jensen's inequalities, together with the fact that $\kappa_1\leq\kappa_2$. Indeed, we note that
\[ 
|\cdot|^{-\gamma}\leq h_1+ h_2,
\]
where $h_1(x)=\begin{cases}|x|^{-\gamma} &\text{if } x\in B(1),\\~0 &\text{if } x\in \RR^d\setminus B(1),  \end{cases}$ and $h_2(x)=\begin{cases}0 &\text{if } x\in B(1),\\1 &\text{if } x\in \RR^d\setminus B(1).  \end{cases}$
Therefore,
\begin{align*}
&\norm{Y(t-s,\cdot)\star\left(|\cdot|^{-\gamma}|u|^{p-1}u(s)\right)}_q\\
&\leq\norm{Y(t-s,\cdot)}_q\norm{h_1|u|^{p-1}u(s)}_1+\norm{Y(t-s,\cdot)}_{q}\norm{h_2|u|^{p-1}u(s)}_{1}\\
&\lesssim (t-s)^{-\frac{\alpha d}{\beta}\left(1-\frac{1}{q}\right)+\alpha-1}\norm{h_1}_{\frac{q}{q-p}}\norm{|u|^{p-1}u(s)}_{\frac{q}{p}}+(t-s)^{-\frac{\alpha d}{\beta}\left(1-\frac{1}{q}\right)+\alpha-1}\norm{|u|^{p-1}(s)}_{\frac{q}{q-1}}\norm{u(s)}_{q}\\
&\lesssim (t-s)^{-\frac{\alpha d}{\beta}\left(1-\frac{1}{q}\right)+\alpha-1}\left(\int_{B(1)}|y|^{-\frac{\gamma q}{q-p}}dy\right)^{\frac{q-p}{q}}\norm{u(s)}_q^p\\
&~~~+(t-s)^{-\frac{\alpha d}{\beta}\left(1-\frac{1}{q}\right)+\alpha-1}\norm{u(s)}_{q}^{p-1}\norm{u(s)}_{q}\\
&= C(\alpha,\beta,d,p,q,\gamma)(t-s)^{-\frac{\alpha d}{\beta}\left(1-\frac{1}{q}\right)+\alpha-1}\norm{u(s)}_q^p\\
&\leq C(\alpha,\beta,d,p,q,\gamma)\norm{u}_E^p(t-s)^{-\frac{\alpha d}{\beta}\left(1-\frac{1}{q}\right)+\alpha-1}s^{-\frac{\alpha d p}{\beta}\left(\frac{1}{q_c}-\frac{1}{q}\right)}.
\end{align*}
By using hypothesis $q<\frac{\alpha dp(p-1)}{\beta(\alpha p -p+1)}$, it follows that
\begin{align*}
J_1(t)&\lesssim\int_{0}^{\frac{t}{2}}(t-s)^{-\frac{\alpha d}{\beta}\left(1-\frac{1}{q}\right)+\alpha-1}s^{-\frac{\alpha d p}{\beta}\left(\frac{1}{q_c}-\frac{1}{q}\right)}ds\\
&\lesssim\left(\frac{t}{2}\right)^{-\frac{\alpha d}{\beta}\left(1-\frac{1}{q}\right)+\alpha-1}
\int_{0}^{\frac{t}{2}}s^{-\frac{\alpha d p}{\beta}\left(\frac{1}{q_c}-\frac{1}{q}\right)}ds\\
&=C_1t^{-\frac{\alpha d}{\beta}\left(1-\frac{1}{q}\right)+\alpha-\frac{\alpha d p}{\beta}\left(\frac{1}{q_c}-\frac{1}{q}\right)}.
\end{align*}
On the other hand, hypothesis $q<\kappa_1$ shows that
\begin{align*}
J_2(t)&\lesssim\int_{\frac{t}{2}}^{t}(t-s)^{-\frac{\alpha d}{\beta}\left(1-\frac{1}{q}\right)+\alpha-1}s^{-\frac{\alpha d p}{\beta}\left(\frac{1}{q_c}-\frac{1}{q}\right)}ds\\
&\lesssim\left(\frac{t}{2}\right)^{-\frac{\alpha d p}{\beta}\left(\frac{1}{q_c}-\frac{1}{q}\right)}
\int_{\frac{t}{2}}^t(t-s)^{-\frac{\alpha d}{\beta}\left(1-\frac{1}{q}\right)+\alpha-1}ds\\
&=C_2t^{-\frac{\alpha d}{\beta}\left(1-\frac{1}{q}\right)+\alpha-\frac{\alpha d p}{\beta}\left(\frac{1}{q_c}-\frac{1}{q}\right)}.
\end{align*}
For $J_3$ we have
$$J_3(t)\lesssim |B|t^{-\frac{\alpha d}{\beta}\left(1-\frac{1}{q}\right)+\alpha-1}$$
and hypothesis $q<\frac{\alpha dp(p-1)}{\beta(\alpha p -p+1)}$ yields 
$$-1<-\frac{\alpha d p}{\beta}\left(\frac{1}{q_c}-\frac{1}{q}\right).$$
Consequently,
\begin{equation}
\label{decayNonLinear}
\norm{\int_{0}^{t} R(t-s)\left(|\cdot|^{-\gamma}|u|^{p-1}u(s)\right) ds-BY(t,\cdot)}_q\lesssim t^{-\frac{\alpha d}{\beta}\left(1-\frac{1}{q}\right)+\alpha-\frac{\alpha d p}{\beta}\left(\frac{1}{q_c}-\frac{1}{q}\right)}
\end{equation}
for $t\geq 1$. This estimate is optimal in the sense that 
$$t^{-\frac{\alpha d}{\beta}\left(1-\frac{1}{q}\right)+\alpha-\frac{\alpha d p}{\beta}\left(\frac{1}{q_c}-\frac{1}{q}\right)}<t^{-\frac{\alpha d}{\beta}\left(\frac{1}{q_c}-\frac{1}{q}\right)},\quad t>1,$$
which is ensured by the assumption $\frac{q}{q-p}<\frac{d}{\gamma}$. Estimates \eqref{decayinit}-\eqref{decayNonLinear} complete the proof.
\end{proof}

\section{Non-existence of positive solutions}
\label{Sec:6}
In this section, we stablish conditions for the non-existence of positive local solutions to the problem \eqref{FE:0} in the sense of Definition \ref{DefSolution}. 
\begin{theorem}
\label{Theorem_blow-up}
Let $\alpha\in (0,1)$, $\beta\in (0,2)$, $0<\gamma <d$ and $q\in [1,\infty)$. Suppose that assumption $(H2)$ holds and that $p > q\left(1+\frac{\beta}{\alpha (d-\gamma)}\right)$. Then there exists a positive function $u_0\in L_q(\mathbb{R}^d)$ such that \eqref{FE:0} possesses no local positive solution. Moreover, any fixed point $u$ of the operator \eqref{OperatorM} is not in $L_{1,\text{loc}}(\mathbb{R}^d)$ for any $t > 0$.
\end{theorem}
\begin{proof}
Let $R>1$ and $\eta >0$ satisfying $0<\eta q<d-\gamma$. We define the function  $$u_0(x)=\begin{cases}|x|^{-\eta} &\text{if } x\in B(R),\\~0 &\text{if } x\in \RR^d\setminus B(R).  \end{cases}$$
We see that  $u_0\in L_1(\RR^d)\cap L_q(\RR^d)$. Suppose that there exists $u\in C([0,T];L_q(\RR^d))$ such that $u(t,x)\geq 0$ a.e. in $(0,T]\times\RR^d$ and that
$$u(t)= S(t)u_0 + \int_{0}^{t} R(t-s)\left(|\cdot|^{-\gamma}u^{p}(s)\right) ds.$$
Properties of operators $S,R$ imply that
$$u(t)\geq S(t)u_0$$
and therefore
$$u(t)\geq \int_{0}^{t} R(t-s)\left(|\cdot|^{-\gamma}(S(s)u_0)^{p}\right) ds. $$
Construction of $f$ in \eqref{phi_i}-\eqref{fOsgood} implies that $(\cdot)^p>Cf$ on $[0,\infty)$ for some positive constant $C$, which depends on $p$ and the chosen sequence $(\phi_i)_{i\in \NN}$. Suppose that $0<t\leq 1$ and take $i$ sufficiently large so that $\left(\frac{\phi_i}{M}\right)^{-\frac{1}{\eta\rho}}\leq t$, with $M, \rho$ as in Proposition \ref{propo:z:phi}. Together with definition of operator $R$ given in \eqref{operatosRS}, we obtain
\begin{align*}
	\int_{B(\epsilon)}u(t,x)dx & \geq \int_{B(\epsilon)}\int_{0}^{t} \int_{\RR^d} Y(t-s,x-y)|y|^{-\gamma}(z(s,y))^p dy\,ds\,dx\\
    & \geq C \int_{0}^{(\phi_i/M)^{-\frac{1}{\eta \rho}}}  \int_{s^{\alpha/\beta}\leq |y|\leq s^{\rho}}|y|^{-\gamma} \int_{B(\epsilon)} Y(t-s,x-y) f(z(s,y))dx\,dy\,ds\\	
	& \geq C \int_{0}^{(\phi_i/M)^{-\frac{1}{\eta \rho}}}  \int_{s^{\alpha/\beta}\leq |y|\leq s^{\rho}} |y|^{-\gamma}\int_{B(\epsilon)} Y(t-s,x-y)\frac{\phi_{i+1}}{2} dx\,dy\,ds\\
	& = C\frac{\phi_{i}^p}{2}\int_{0}^{(\phi_i/M)^{-\frac{1}{\eta \rho}}}  \int_{s^{\alpha/\beta}\leq |y|\leq s^{\rho}} |y|^{-\gamma}\int_{|x+y|\leq\epsilon} Y(t-s,x)dx\,dy\,ds.
\end{align*}
From now on, the rest of the proof includes argument similar to those used in \cite[Theorem 3.1]{SoVe25}, the only difference being the factor $|\cdot|^{-\gamma}$. In this way,
\begin{align*}
	\int_{B(\epsilon)}u(t,x)dx & \geq C(d,\epsilon,\alpha,p) t^{\alpha-1}\left(\phi_{i}^{p - \frac{1}{\eta \rho}((d-\gamma)\rho + 1)}-\phi_{i}^{p - \frac{1}{\eta \rho}((d-\gamma)\alpha/\beta + 1)}\right).
\end{align*}
Thereby, if $p>q\left(1+\frac{\beta}{\alpha (d-\gamma)}\right)$ then it is possible to choose $\eta \in\left(0,\frac{d-\gamma}{q}\right)$ and $\rho\in\left(0,\frac{\alpha}{\beta}\right)$ such that $p > \frac{1}{\eta\rho}\left((d-\gamma)\rho + 1\right)$. Consequently, 
\[
\int_{B(\epsilon)}u(t,x)dx \rightarrow \infty
\]
whenever $i\rightarrow \infty$. This shows that the blow-up is local to the origin which leads to a contradiction.
\end{proof}
Next, for $R>1$ we define the function
$$\varphi(x)=\begin{cases}|x|^{-\frac{d}{q_c}} &\text{if } x\in B(R),\\~0 &\text{if } x\in \RR^d\setminus B(R),  \end{cases}$$
where $q_c$ is the critical exponent given in \eqref{Cexponent}. By employing this function, the following result put forward a situation contrary to that presented in Corollary \ref{corol_global}. 
\begin{corollary}
\label{corol_BlowUp}
Let $\alpha\in (0,1)$, $\beta\in (0,2)$, $0<\gamma<\min(\beta,d)$ and $q\in [1,\infty)$. Suppose that assumption $(H2)$ holds and that $p > q\left(1+\frac{\beta}{\alpha (d-\gamma)}\right)$. If $u_0\geq B\varphi$ for some  constant $B>0$, then the problem \eqref{FE:0} possesses no local positive solution. Moreover, any fixed point $u$ of the operator \eqref{OperatorM} is not in $L_{1,\text{loc}}(\mathbb{R}^d)$ for any $t > 0$.
\end{corollary}
\begin{proof}
Hypotheses on $q$ and $p$ show that  $$\frac{\beta-\gamma}{p-1}q<d-\gamma$$ and the conclusion follows from Theorem \ref{Theorem_blow-up} taking $\eta=\frac{\beta-\gamma}{p-1}=\frac{d}{q_c}$.
\end{proof}

\section{Discussion}
\label{Sec:7}
Abundant works have been developed on the classical Hardy parabolic equation \eqref{FE:0:1:2}. However, in this section we limit our discussion to those studied in \cite{STW17} and \cite{HiSi24}. Our results for global solutions, stated in Theorems \ref{globalsolutionsmallu0} and \ref{Global2}, provide a natural extension of \cite[Theorem 1.3]{STW17} for the solvability of \eqref{FE:0} in $L_q(\RR^d)$ spaces. Under $|u_0(x)|\leq A|x|^{-\frac{\beta-\gamma}{p-1}}$ and a small constant $A$, Theorem \ref{Global2} is a similar result to that found in \cite[Theorem 1.3 (iii)]{STW17} with $\alpha=1$, $\beta=2$ and $\omega_{\nu}\equiv 1$, if $p>1+\frac{\beta-\gamma}{d}$ equivalent to $q_c>1$. Theorem \ref{globalsolutionsmallu0} also generalizes \cite[Theorem 1.3 (i)]{STW17}. On the other hand, Corollary \ref{corol_BlowUp} claims that problem \eqref{FE:0} possesses no local-in-time positive solutions whenever $u_0(x)\geq B|x|^{-\frac{\beta-\gamma}{p-1}}$ for a suitable constant $B>0$, the same being true in \cite{HiSi24} with $\alpha=1$, $\beta\in (0,2)$ and $\omega_{\nu}\equiv 1$, although in our case with a different condition for $p$ but which leads to $p>1+\frac{\beta-\gamma}{d}$. 

The study of solutions to \eqref{FE:0} in other spaces, such as $C_0(\RR^d)$ (the space of continuous functions on $\RR^d$ vanishing at infinity), remains open. Under assumption (H2), it is a known fact that the operator $-\Psi_{\beta}(-i\nabla)$ is the generator of a Feller semigroup, which could help to develop new methods and improve some results of this work.

\bibliographystyle{plain}
\bibliography{RefHardy-typeEq}

\end{document}